\documentclass[11pt]{article}
\makeatletter \oddsidemargin .4in \evensidemargin -.1in \textwidth
14cm \topmargin 0.1cm \textheight 22cm
\newcommand{\singlespacing}{\let\CS=\@currsize\renewcommand{\baselinestretch}{1}\tiny\CS}
\newcommand{\doublespacing}{\let\CS=\@currsize\renewcommand{\baselinestretch}{1.35}\tiny\CS}
\newcommand{\doublespacinggg}{\let\CS=\@currsize\renewcommand{\baselinestretch}{1.55}\tiny\CS}
\newcommand{\doublespacingg}{\let\CS=\@currsize\renewcommand{\baselinestretch}{1.75}\tiny\CS}
\doublespacinggg

\newtheorem{thm}{Theorem}[section]
 \newtheorem{cor}[thm]{Corollary}
 \newtheorem{lem}[thm]{Lemma}

\newcommand{\be}{\begin{equation}}
\newcommand{\ee}{\end{equation}}
\newcommand{\bea}{\begin{eqnarray}}
\newcommand{\eea}{\end{eqnarray}}
\newcommand{\nn}{\nonumber}
\newcommand{\bee}{\begin{eqnarray*}}
\newcommand{\eee}{\end{eqnarray*}}
\newcommand{\lb}{\label}
\newcommand{\nii}{\noindent}

\title{ Generalized derivations as a generalization of Jordan homomorphisms acting on Lie ideals and right ideals}
\author{\small   Basudeb Dhara, Shervin Sahebi and Venus Rahmani}

\begin{document}
\date{}
\maketitle \noindent \vspace{-.8cm}

\doublespacingg

\begin{center}
\begin{minipage}{11cm} \footnotesize { \textsc{Abstract:}
Let $R$ be a prime ring with center $Z(R)$ and extended centroid
$C$, $H$ a non-zero generalized derivation of $R$ and $n\geq 1$ a
fixed integer. In this paper we study the situations: (1)
$H(u^2)^n-H(u)^{2n}\in C$ for all $u \in L$, where $L$ is a
non-central Lie ideal of $R$; (2) $H(u^2)^n-H(u)^{2n}=0$ for all
$u\in [I, I]$, where $I$ is a nonzero right ideal of $R$.}

 \end{minipage}
\end{center}

 \vspace*{.4cm}

 \noindent {\footnotesize {\bf Mathematics Subject Classification 2010:} 16N60, 16U80, 16W25.  }

\vspace*{.2cm}

\doublespacing

\section{Introduction}

Throughout this paper, $R$ always denotes a prime ring with center
$Z(R)$ and with extended centroid $C$, $U$ the Utumi quotient ring
of $R$. For given $x,y\in R$, the Lie commutator of $x,y$ is
denoted by $[x,y]$ and defined by $[x,y]=xy-yx$. A linear mapping
$d: R\rightarrow R$ is called a derivation, if it satisfies the
Leibniz rule $d(xy)=d(x)y+ xd(y)$ for all $x, y \in R$.
  We recall that an additive map $H: R\rightarrow R$ is called a generalized
 derivation, if there exists a derivation $d: R\rightarrow R$ such
 that $H(xy)=H(x)y+xd(y)$ holds for all $x, y \in R$.
 Let $S$ be a nonempty subset of $R$ and $F:R\rightarrow R$ be an
additive mapping. Then we say that $F$ acts as homomorphism or
anti-homomorphism on $S$ if $F(xy)=F(x)F(y)$ or $F(xy)=F(y)F(x)$
holds for all $x,y\in S$ respectively. The additive mapping $F$
acts as a Jordan homomorphism on $S$ if $F(x^2)=F(x)^2$ holds for
all $x\in S$.

 Several authors studied the situations, when some specific type
 of additive maps acts as homomorphisms or anti-homomorphisms in
 some subsets of $R$.
For instance Asma, Rehman and Shakir in~\cite{a01} proved that if
$d$ is a derivation of a  $2$-torsion free prime ring $R$ which
acts as a homomorphism or ani-homomorphism on a square closed Lie
ideal $L$ of $R$, then $d=0$ or $L\subseteq Z(R)$. Recently,
in~\cite{Gol} Golbasi and Kaya study the case when derivation $d$
is replaced by generalized derivation $H$. More precisely, they
proved the following: Let $R$ be a prime ring of characteristic
different from $2$, $H$ a generalized derivation of $R$, $L$ a Lie
ideal of $R$ such that $u^2\in L$ for all $u\in L$. If $H$ acts as
a homomorphism or anti-homomorphism on $L$, then either $d=0$ or
 $L\subseteq Z(R)$.

Recently in \cite{filippis}, De Filippis studied the situation
when generalized derivation $H$ acts as a Jordan homomorphism on a
non-central Lie ideal $L$ of $R$ and on the set $[I,I]$, where $I$
is a nonzero right ideal of a prime ring $R$.

In the present paper our motivation is to generalize all the above
results by studying the following situations: (1)
$H(u^2)^n-H(u)^{2n}\in C$ for all $u \in L$, where $L$ is a
non-central Lie ideal of $R$; (2) $H(u^2)^n-H(u)^{2n}=0$ for all
$u\in [I, I]$, where $I$ is a nonzero right ideal of $R$.

 \noindent
The following results are useful tools
 needed in the proof of main results.

\noindent

 {\bf Remark 1.} Let $R$ be a prime ring and $L$ a noncentral
Lie ideal of $R$. If char$(R)\neq 2$, by \cite[Lemma 1]{a010}
there exists a nonzero ideal $I$ of $R$ such that $0 \neq [I,R]
\subseteq L$. If char$(R)=2$ and dim$_CRC>4$, i.e., char$(R)=2$
and $R$ does not satisfy $s_4$, then by \cite[Theorem 13]{lm}
there exists a nonzero ideal $I$ of $R$ such that $0 \neq [I,R]
\subseteq L$. Thus if either char$(R)\neq 2$ or $R$ does not
satisfy $s_4$, then we may conclude that there exists a nonzero
ideal $I$ of $R$ such that $[I,I] \subseteq L$.

\noindent

 {\bf Remark 2.} Let $R$ be a prime ring and $U$ be the Utumi quotient
ring of $R$ and $C=Z(U)$, the center of $U$.
It is well known that any derivation of $R$ can be uniquely
extended to a derivation of $U$, In~\cite{g}
Lee proved that every generalized derivation $H$ on a dense
right ideal of $R$ can be uniquely extended to a generalized derivation
of $U$ and assume the form $H(x)=ax+d(x)$ for all $x\in U$
 ,some $a\in U$ and a derivation $d$ of $U$.

%%%%%%%%%--------------------LEMMA-----------------------------------
\section{Generalized derivations on Lie ideals}

 \noindent
 We establish the following results required in the proof of Theorem~\ref{THM1}.

 \begin{lem}\label{lem1.1}
 Let $R=M_k(F)$, be the ring of all $k\times k$ matrices
 over a field $F$ with $k\geq 2$, $a\in R$ and $n\geq 1$ a fixed integer.
 If $(a[x,y]^2)^n-(a[x,y])^{2n}=0$ for all $x,y \in R$, then $a\in F\cdot I_k$ and either  $a=0$ or $a^n=1$.
 \end{lem}
 %%%%%%--------------PROOF---------------------------------------
 {\it Proof.}
 Let $a=(a_{ij})_{k\times k}$ where $a_{ij}\in F$.
  By choosing $x=e_{ii}$, $y=e_{ij}$ for any $i\neq j$, we have
 \begin{equation}\label{b1}
 0=-(ae_{ij})^{2n}.
 \end{equation}
  Left multiplying (1) by $e_{ij}$,
 it gives $$0=e_{ij}(ae_{ij})^{2n}=a_{ji}^{2n}e_{ij},$$
 implying $a_{ji}=0$.  Thus for any $i\neq j$, we have $a_{ij}=0$,
 which implies that $a$ is a diagonal matrix. Let $a=\sum_{i=1}^k a_{ii}e_{ii}$.
 For any $F$-automorphism
 $\theta$ of $R$, we have $(a^{\theta}[x,y]^2)^n-(a^{\theta}[x,y])^{2n}=0$
 for every $x, y \in R$. Hence $a^{\theta}$
 must also be diagonal. We have
 $$(1+e_{ij})a(1-e_{ij})=\sum_{i=1}^k a_{ii}e_{ii}+(a_{jj}-a_{ii})e_{ij}$$
 diagonal. Therefore, $a_{jj}=a_{ii}$ and so $a\in F\cdot I_k$.
Thus the main assumption reduces to
$$a^n(a^{n}-1)[x,y]^{2n}=0$$ for all $x,y \in R$. By choosing
$x=e_{ij}, y=e_{ji}$ we get
$0=a^n(a^{n}-1)[e_{ij},e_{ji}]^{2n}=a^n(a^{n}-1)\{e_{ii}+e_{jj}\}.$
This leads either $a=0$ or $a^n=1$.

%%%%%%%%%%%%%%%%%%%%%%%%         LEMMA 2          %%%%%%%%%%%%%%%%%%%%%%%%%%%%%%%%%%%%%%%%%%%%%%%%

\begin{lem}\label{b1}
 Let $R=M_k(F)$ be the ring of all $k\times k$ matrices
 over a field $F$ with $k\geq 3$, $a, b\in R$ and $n\geq 1$ a fixed integer.
 If $(a[x,y]^2-[x,y]^2b)^n-(a[x,y]-[x,y]b)^{2n}\in F\cdot I_k$,
 for all $x,y \in R$, then $a,b\in F\cdot I_k$ and $a-b=0$ or $(a-b)^n=1$.
 \end{lem}
 %%%%%%--------------PROOF---------------------------------------
 {\it Proof.}
 Let $a=(a_{ij})_{k\times k}$ and $b=(b_{ij})_{k\times k}$ where $a_{ij}, b_{ij}\in F$.
By assumption we have
$$[(a[x,y]^2-[x,y]^2b)^n-(a[x,y]-[x,y]b)^{2n},z]=0,$$
for all $x,y,z\in R$.
 By choosing
 $x=e_{ii}$, $y=e_{ij}$ and $z=e_{ik}$ for any
 $i\neq j\neq k$, we have
$$0=[(ae_{ij}-e_{ij}b)^{2n},e_{ik}]=(e_{ij}b)^{2n}e_{ik}-e_{ik}(ae_{ij})^{2n}=(b_{ji})^ne_{ik}-a_{ki}(a_{ji})^{2n-1}e_{ij}.$$
 Thus $b_{ji}=0$.  We conclude that $b$ is a diagonal matrix.
 By the same argument in Lemma \ref{lem1.1}, we have $b\in F\cdot I_k$.
 Similarly we can conclude $a\in F\cdot I_k$.
 Therefore the main assumption says that
 $$(a-b)^n(1-(a-b)^n)([[x,y]^{2n},z])=0.$$
 Hence $a-b=0$ or $(a-b)^n=1$.

%%%%%%%%%%%%%%%%%%%%%%%%         LEMMA 3          %%%%%%%%%%%%%%%%%%%%%%%%%%%%%%%%%%%%%%%%%%%%%%%%
  \begin{lem}\label{b2}
 Let $R$ be a noncommutative prime ring with extended centroid
 $C$, $I$ a nonzero ideal of $R$ and $a, b\in R$.
 Suppose that  $(a[x,y]^2-[x,y]^2b)^n=(a[x,y]-[x,y]b)^{2n}$
 for all $x,y \in I$, where $n\geq 1$ is a fixed integer.
  Then $a, b\in C$ and either $a-b=0$ or $(a-b)^n=1$.
 \end{lem}
 %%%%%%%%%%%%%%%%%%-------PROOF-----------------------------------------
 {\it Proof.}
By assumption, $I$ satisfies the generalized polynomial identity
$$F(x,y)=(a[x,y]^2-[x,y]^2b)^n-(a[x,y]-[x,y]b)^{2n}.$$
By Chuang \cite[Theorem 2]{a1}, this
generalized polynomial identity (GPI) is also satisfied by $U$.
If $a\notin C$ or $b\notin C$, then $F(x,y)=0$ is a nontrivial
 (GPI) for $U$. In case $C$ is infinite, we have $F(x,y)=0$ for
 all $x, y \in U\bigotimes_C\overline{C}$
where $\overline{C}$ is the algebraic closure of $C$.
 Since both $U$ and $U\bigotimes_C\overline{C}$
are prime and centrally closed \cite{u54}, we may replace $R$ by
$U$ or $U\bigotimes_C\overline{C}$ according to $C$ is finite or
infinite. Thus we may assume that $R$ is centrally closed over $C$
which is either finite or algebraically closed and $F(x,y)=0$ for
all $x, y\in R$. By Martindale's Theorem \cite{a6}, $R$ is then a
primitive ring having nonzero $soc(R)$ with $C$ as the
associated division ring. Hence by Jacobson's Theorem \cite{a5},
$R$ is isomorphic to a dense ring of linear transformations of a
vector space $V$ over $C$. If dim$_CV=k$, then the density of $R$
on $V$ implies that $R\cong M_k(C)$. Since $R$ is noncommutative,
$k\geq 2$.

 We want to show that for any $v\in V$,
  $v$ and $bv$ are linearly $C$-dependent.
  Suppose on contrary that $v$ and $bv$ are linearly $C$-independent for some
  $v\in V$. By density there exist $x, y \in R$ such that
  $$\begin{array}{ccc}
    xv=0,  & xbv=-bv, \\
    yv=v,  & ybv=v.
  \end{array}$$
  Then $[x,y]v=0$, $[x,y]bv=v$, and so  $[x,y]^2bv=0$.
  Hence
  $$0=((a[x,y]^2-[x,y]^2b)^n-(a[x,y]-[x,y]b)^{2n})v=-v,$$
  a contradiction.
  Thus we conclude that $\{v, bv\}$ is a linearly $C$-dependent set of vectors for any $v\in
  V$. Thus for any $v\in V$, $bv=\alpha_v v$ for some
$\alpha_v\in C$. Now we prove that $\alpha_v$ is independent of
the choice of $v\in V$. Let $u$ be a fixed vector of $V$. Then
$bu=\alpha u$. Let $v$ be any vector of $V$. Then $bv= \alpha_v
v$, where $\alpha_v\in C$.
 If $u$ and $v$ are linearly $C$-dependent, then $u=\beta v$,
for $\beta\in C$. In this case, we see that $\alpha u=bu=\beta bv=
\beta(\alpha_v v) = \alpha_v(\beta v) =\alpha_v u$, implying
$\alpha=\alpha_v$.

Now if $u$ and $v$ are linearly $C$-independent, then we have
$\alpha_{u+v}(u+v)=b(u+v)=bu+bv=\alpha u+\alpha_v v$, which
implies $(\alpha_{u+v}-\alpha)u+(\alpha_{u+v}-\alpha_v)v=0$. Since
$u$ and $v$ are linearly $C$-independent, we  have
$\alpha_{u+v}-\alpha=0=\alpha_{u+v}-\alpha_v$ and so
$\alpha=\alpha_v$. Thus $bv = \alpha v$ for all $v\in V$, where
$\alpha\in C$ is independent of the choice of $v\in V$.

 Now, let $r\in R$ and $v\in V$. Since $bv=\alpha v$,
 $$[b,r]v=(br)v-(rb)v=b(rv)-r(bv)=(rv)\alpha-r(v\alpha)=0,$$
 that is $[b,r]V=0$.
 Hence $[b,r]=0$ for all $r\in R$, implying $b\in C$.

 Then our assumption reduces to
 $(a'[x,y]^2)^n-(a'[x,y])^{2n}=0$ for all $x,y\in R$, where
 $a'=a-b$. If dim$_CV=k$, then by Lemma \ref{lem1.1},
 we have $a'=a-b\in C$ and either $a'=0$ or $a'^n=1$. Since $b\in C$, $a\in
 C$. Let dim$_CV=\infty$. Then
for any $e^2=e\in soc(R)$ we have $eRe\cong M_t(C)$ with
$t=$dim$_CVe$. Assume that $a'\notin C$. Then $a$ does not
centralize the nonzero ideal $soc(R)$. Hence there exist $h\in
soc(R)$ such that $[a,h]\neq 0$. By Litoff's theorem \cite{x11},
there exists idempotent $e\in soc(R)$ such that $a'h, ha',h\in
eRe$. We have $eRe\cong M_k(C)$ with $k=$ dim$_CVe$. Since $R$
satisfies generalized identity
$e\{(a'[exe,eye]^2)^n-(a'[exe,eye])^{2n}\}e=0$,
 the subring $eRe$ satisfies
$(ea'e[x,y]^2)^n-(ea'e[x,y])^{2n}=0$. Then by the above finite
dimensional case, $ea'e$ is a central element of $eRe$. Thus
$ah=(eae)h=heae=ha$, a contradiction. Hence we conclude that
$a'\in C$. Then our identity reduces to $a'^n(a'^n-1)[x,y]^{2n}=0$
for all $x,y\in R$. Since dim$_CV=\infty$, $R$ can not satisfy any
polynomial identity, and hence $a'^n(a'^n-1)=0$ implying either
$a'=0$ or $a'^n=1$. Since $a'=a-b$, we obtain our conclusion.

 %%%%%%%----------------------------------------------------------

 \begin{thm}\label{THM1}
 Let $R$ be a prime ring, $H$ a nonzero
 generalized derivation of $R$ and $L$ a non-central Lie ideal of $R$.
 Suppose that $H(u^2)^n-H(u)^{2n}=0$ for all
 $u \in L$, where $n\geq 1$ is a fixed integer.
 Then one of the following holds:
  \begin{enumerate}
   \item $char (R)=2$ and $R$ satisfies $s_4$;
   \item  $H(x)=bx$ for some $b\in C$ and $b^n=1$.
 \end{enumerate}
 \end{thm}

 {\it Proof.}
 We assume that either char$(R)\neq 2$ or $R$ does not satisfy $s_4$.
 Since $L$ is non central by Remark 1,
 there exists a nonzero ideal $I$ of $R$ such that $[I,I]\subseteq L$.
  Thus by assumption, $I$ satisfies the differential identity
 $$H([x,y]^2)^n-H([x,y])^{2n}=0.$$
 Since $I$ and $U$ satisfy the same differential identities \cite{g}, we may assume that
 $H([x,y]^2)^n-H([x,y])^{2n}=0$ for all $x, y \in U$.
 As we have remarked in Remark 2,
 we may assume that for all $x\in U$, $H(x)=bx+d(x)$ for some $b\in U$ and a derivation
 $d$ of $U$. Hence $U$ satisfies
\begin{equation}\label{b4}
(b[x,y]^2+d([x,y]^2])^n-(b[x,y]+d([x,y]])^{2n}=0.
 \end{equation}
 Assume first that $d$ is inner derivation of $U$, i.e., there
 exists $p\in U$ such that $d(x)=[p,x]$ for all $x\in U$. Then
 $$(b[x,y]^2+[p,[x,y]^2])^n-(b[x,y]+[p,[x,y]])^{2n}=0,$$
for all $x, y\in U$ that is
$$((b+p)[x,y]^2-[x,y]^2p)^n-((b+p)[x,y]-[x,y]p)^{2n}=0,$$
for all $x, y\in U$. By Lemma~\ref{b2}, $b+p, p\in C$ and $b=0$ or
$b^n=1$. If $b=0$ then $H(x)=0$, a contradiction. Otherwise,
$H(x)= bx$ for some $b\in C$ and $b^n=1$, as  desired.

\noindent
On the other hand (2) implies
$$\begin{array}{lr}
             (b[x,y]^2+([d(x),y]+[x,d(y)])[x,y]+[x,y]([d(x),y]+[x,d(y)]))^n &  \\
             -(b[x,y]+[d(x),y]+[x, d(y)])^{2n}=0,&
           \end{array}$$ for all $x,y\in U$.
So if $d$ is not $U$-inner, then by Kharchenko's
theorem~\cite{a3}, we have
 $$\begin{array}{lr}
 (b[x,y]^2+([z,y]+[x,t])[x,y]+[x,y]([z,y]+[x,t]))^n  & \\
 -(a[x,y]+[z,y]+[x,t])^{2n}=0, &
  \end{array}$$
for all $x, y, z, t\in U$. In particular, for $x=t=0$, we have
$[z,y]^{2n}=0$ for all $z,y \in U$. Note that this is a polynomial
identity and hence there exists a field
 $F$ such that $R\subseteq M_k(F)$, the ring of $k\times k$
 matrices over a field $F$, where $k\geq 1$.
 Moreover, $R$ and $M_k(F)$ satisfy the same polynomial identity~\cite[Lemma 1]{a03}
 that is $[z,y]^{2n}=0$ for all $y,z\in M_k(F)$.
 But by choosing $z=e_{12}$, $y=e_{21}$ we get
 $$0=[z,y]^{2n}=e_{11}+e_{22}$$
 which is a  contradiction.
%%%%%%%%%%%%%%%%%%%%%%%%%%%%%%%%%%%%%%%%%%%%%%%THEOREM2%%%%%%%%%%%%%%%%%%%%%%%%%%%%%%%%%%%%%%%

\begin{lem}\label{lab1}
  Let $R$ be a noncommutative prime ring with extended centroid
 $C$ and $a, b\in R$.
 Suppose that $(a[x,y]^2-[x,y]^2b)^n-(a[x,y]-[x,y]b)^{2n}\in C$
 for all $x,y \in R$, where $n\geq 1$ is a fixed integer.
  Then one of the following holds:
   \begin{enumerate}
  \item  $a, b\in C$, such that $a-b=0$ or $(a-b)^n=1$;
  \item  $R$ satisfies $s_4$.
  \end{enumerate}
\end{lem}
 {\it Proof.} Since $R$ and $U$ satisfy the same generalized
 polynomial identities  (see \cite{a1}), $U$ satisfies
  \bea g(x,y,z)=[(a[x,y]^2-[x,y]^2b)^n-(a[x,y]-[x,y]b)^{2n},z].
  \eea
 Suppose first that $g(x,y,z)$ is a trivial generalized polynomial
 identity for $R$. Let $T=U\ast_CC\{x,y,z\}$ be the free product
 of $U$ and $C\{x,y,z\}$, the free $C$-algebra in noncommuting
 indeterminates $x,y,z$. Then
 $$[(a[x,y]^2-[x,y]^2b)^n-(a[x,y]-[x,y]b)^{2n},z]$$ is zero element
 in $T$. Let $a\notin C$. Then $a$ and $1$ are linearly
 independent over $C$. Thus from above,
 $$\{a[x,y]^2(a[x,y]^2-[x,y]^2b)^{n-1}-a[x,y](a[x,y]-[x,y]b)^{2n-1}\} z$$ is zero element
 in $T$ that is
 $$a[x,y]\bigg\{[x,y](a[x,y]^2-[x,y]^2b)^{n-1}-(a[x,y]-[x,y]b)^{2n-1}\bigg\}
 z=0$$ in $T$. Again since $a$ and $1$ are linearly independent, we
 have $$a[x,y]\bigg\{-a[x,y](a[x,y]-[x,y]b)^{2n-2}\bigg\}
 z=0$$ and so $a[x,y]\{-a[x,y](a[x,y])^{2n-2}\} z=0$ in $T$ implying
 $a=0$, a contradiction. Hence $a\in C$. Then the identity reduces
 to $$[([x,y]^2(a-b))^n-([x,y](a-b))^{2n},z]=0.$$ Again if $a-b\notin C$, then it gives
 $$z\bigg\{([x,y]^2(a-b))^{n-1}[x,y]^2(a-b)-([x,y](a-b))^{2n-1}[x,y](a-b)\bigg\}=0$$
 that is
 $$z\bigg\{([x,y]^2(a-b))^{n-1}[x,y]-([x,y](a-b))^{2n-1}\bigg\}[x,y](a-b)=0$$ in
 $T$. This again implies $z\{-([x,y](a-b))^{2n-1}\}[x,y](a-b)=0$,
 implying $a-b=0$, a contradiction. Hence $a-b\in C$. Since $a\in
 C$, we have $b\in C$. Then the (GPI) becomes $(a-b)^n((a-b)^n-1)[x,y]^{2n}\in
 C$. This gives either $a-b=0$ or $(a-b)^n=1$, which is our
 conclusion.

 Next we assume that $g(x,y,z)$ is a nontrivial generalized
 polynomial identity for $R$ and so for $U$. Let $I$ be a
 two-sided ideal of $U$. If
 $(a[x,y]^2-[x,y]^2b)^n-(a[x,y]-[x,y]b)^{2n}=0$ for all $x,y\in
 I$, then the conclusion follows by Lemma~\ref{b2}. Hence we assume that
 there exist $x,y\in I$, such that
 $0\neq (a[x,y]^2-[x,y]^2b)^n-(a[x,y]-[x,y]b)^{2n}\in I\cap C$.
 Then by \cite[Theorem 1]{cl}, $R$ is a PI-ring, therefore $RC=Q=U$ is
 a is a finite-dimensional central simple $C$-algebra by Posner's
  theorem for prime PI-ring. Then by Lemma 2 in \cite{a03}, there exists a field $F$ such that
 $U\subseteq M_k(F)$, the ring of all $k\times k$ matrices over
 $F$, moreover $U$ and $M_k(F)$ satisfy the same generalized
 identities. Therefore $M_k(F)$ satisfies $g(x,y,z)$ and then the
 result follows from Lemma \ref{b1}.

%%%%%%%%%%%%%%%%%%%%%%%%%%%%%%%%%THEOREM 2 %%%%%%%%%%%%%%%%%%%%%%%%
 \noindent
  Now  we are ready to prove Theorem~\ref{THM2}.
 \begin{thm}\label{THM2}
 Let $R$ be a prime ring with extended centroid $C$,
 $H$ a nonzero generalized derivation of $R$ and $L$ a non-central Lie ideal of $R$.
 Suppose that $H(u^2)^n-H(u)^{2n}\in C$ for all
 $u \in L$, where $n\geq 1$ is a fixed integer.
 Then $R$ satisfies $s_4$ or $H(x)=bx$ for some $b\in C$ and $b^n=1$.
 \end{thm}

 {\it Proof.}
 Let $R$ does not satisfy $s_4$. Then by Remark 1,
 there exists an ideal $0\neq I$ of $R$ such that $0\neq [I,I]\subseteq
 L$. Then by assumption, $H([x,y]^2)^n-H([x,y])^{2n}\in C$ for all
 $x,y \in I$. If $H$ is inner generalized derivation of $R$, then
 the result follows by Lemma \ref{lab1}. Let $H$ be not inner. Then by
 Remark 2, $H$ has the form $H(x)=bx+d(x)$, where $b\in U$ and $d$ is
 a derivation of $U$. Since $I$ and $U$ satisfy the same
 generalized polynomial identities (see \cite{a1}) as well as the
 same differential identities (see \cite{g}), we may assume that
 $U$ satisfies
 $[(b[x,y]^2+d([x,y]^2])^n-(b[x,y]+d([x,y]])^{2n},w]=0.$
 Since $H$ is not inner, $d$ is also not inner derivation of $U$.
 We have
$$\begin{array}{lr}
            [(b[x,y]^2+([d(x),y]+[x,d(y)])[x,y]+[x,y]([d(x),y]+[x,d(y)]))^n &  \\
             -(b[x,y]+[d(x),y]+[x, d(y)])^{2n},w]=0.&
           \end{array}$$
By Kharchenko's theorem~\cite{a3} and then by same argument of
Theorem \ref{THM1}, we have
 $[[z,y]^{2n},w]=0$ for all $z,y,w \in U$. This is a polynomial
 identity for $U$. Then by \cite[Lemma 2]{a03}, there exists a field $F$ such that
 $U\subseteq M_k(F)$, the ring of all $k\times k$ matrices over
 $F$, moreover $U$ and $M_k(F)$ satisfy the same generalized
 identities. If $k\leq 2$, then $U$ and so $R$ satisfies $s_4$, as desired. If
 $k\geq 3$, then
 $0=[[z,y]^{2n},w]=[[e_{12},e_{21}]^{2n},e_{13}]=e_{13}$, a
 contradiction.

%%%%%%%%%%%%%%%%%%%%%%%%%%%%%%%%%%SECTION%%%%%%%%%%%%%%%%%%%%%%%%%%%%%%%%%%%%%

\section{Generalized derivations on right ideals}
In this section we will prove the following theorem:
\begin{thm}\label{right}
Let $R$ be a prime ring, $I$ a non-zero right ideal of $R$ and $H$
a non-zero generalized derivation of $R$. If
$H(u^2)^n-H(u)^{2n}=0$ for all $u\in [I, I]$ then one of the
following holds:
\begin{enumerate}
  \item $[I, I]I=0$;
  \item there exists $a\in U$ such that $H(x)=xa$ for all $x\in I$ with $aI=0$;
  \item there exists $a\in U$ such that $H(x)=ax$ for all $x\in R$ with $aI=0$;
  \item there exists $a,b\in U$ such that $H(x)=ax+xb$ for all $x\in R$ with $(a-\alpha)I=(b-\beta)I=0$ for some $\alpha, \beta \in
  C$ and $(\alpha+\beta)^n=1$.
\end{enumerate}
\end{thm}

To prove this theorem, we need the following:
%%%%%%%%%%%%%%%%%%%%%%%%%%%%%%%%LEMMA%%%%%%%%%%%%%%%%%%%%%%%%%%%%%%%%%%%%%%%%
 \begin{lem}\label{c1}
 Let $R$ be a prime ring with extended centroid $C$
 and $I$ a nonzero right ideal of $R$. If for some $a, b\in R$,
 $(a[x_1, x_2]^2+[x_1, x_2]^2b)^n-(a[x_1, x_2]+[x_1,
x_2]b)^{2n}=0$ for all $x_1,x_2\in I$, then $R$ satisfy a
non-trivial generalized polynomial identity or there exist
$\alpha, \beta\in C$ such that $(a-\alpha)I=0$, $(b-\beta)I=0$
with $\alpha+\beta=0$ or $(\alpha+\beta)^n=1$ or $b=-\alpha\in C$.
\end{lem}
%%%%%%%%%%%%%%%%%%%%%%%%%%%PROOF%%%%%%%%%%%%%%%%%%%%%%%%%%%%%%%%%%%%%%%%%%%%
{\it Proof.} By our hypothesis, for any $x_0\in I$, $R$ satisfies
the following generalized identity \bea\lb{hf} (a[x_0x_1,
x_0x_2]^2+[x_0x_1, x_0x_2]^2b)^n-(a[x_0x_1, x_0x_2]+[x_0x_1,
x_0x_2]b)^{2n}.\eea We assume that this is a trivial (GPI) for $R$,
for otherwise we are done. If there exists $x_0\in I$ such that
$\{x_0, ax_0\}$ is linearly $C$-independent, then from above we
have that $R$ satisfies
\bea a[x_0x_1, x_0x_2]^2(a[x_0x_1, x_0x_2]^2+[x_0x_1, x_0x_2]^2b)^{n-1}\nn\\
-a[x_0x_1, x_0x_2](a[x_0x_1, x_0x_2]+[x_0x_1,
x_0x_2]b)^{2n-1},\eea that is
 \bea a[x_0x_1, x_0x_2]\bigg\{[x_0x_1, x_0x_2](a[x_0x_1, x_0x_2]^2+[x_0x_1, x_0x_2]^2b)^{n-1}\nn\\
 -(a[x_0x_1, x_0x_2]+[x_0x_1,x_0x_2]b)^{2n-1}\bigg\}.\eea
 Again since $\{x_0, ax_0\}$ is linearly $C$-independent we have
  $$a[x_0x_1, x_0x_2]\bigg\{-a[x_0x_1, x_0x_2](a[x_0x_1, x_0x_2]+[x_0x_1,
x_0x_2]b)^{2n-2}\bigg\}=0$$ and then by the same manner we have
 $$a[x_0x_1, x_0x_2]\bigg\{-a[x_0x_1, x_0x_2](a[x_0x_1,
 x_0x_2])^{2n-2}\bigg\}=0,$$ which is nontrivial, a contradiction. Thus
 $\{x, ax\}$ is linearly $C$-dependent for all $x\in I$ that is $(a-\alpha)I=0$ for some $\alpha\in C$.
 Then our generalized identity reduces to
$$(\alpha[x_0x_1, x_0x_2]^2+[x_0x_1, x_0x_2]^2b)^n-(\alpha[x_0x_1, x_0x_2]+[x_0x_1,
x_0x_2]b)^{2n}=0$$ that is
 \bea\lb{ps1} ([x_0x_1,
 x_0x_2]^2(b+\alpha))^n-([x_0x_1,x_0x_2](b+\alpha))^{2n}=0.\eea This
 is $$[x_0x_1, x_0x_2]\bigg\{[x_0x_1, x_0x_2](b+\alpha)([x_0x_1, x_0x_2](b+\alpha))^{n-1}-((b+\alpha)[x_0x_1,x_0x_2])^{2n-1}(b+\alpha)\bigg\}=0.$$
If $\{x_0, (b+\alpha)x_0\}$ is linearly independent over $C$, then
 $$[x_0x_1,
 x_0x_2]\bigg\{-((b+\alpha)[x_0x_1,x_0x_2])^{2n-1}(b+\alpha)\bigg\}=0,$$
 which is nontrivial, a contradiction. Thus $\{x, (b+\alpha)x\}$ is linearly dependent over
 $C$ for all $x\in I$, that is $(b+\alpha-\gamma)I=0$ for some $\gamma\in
 C$. Let $\beta=\gamma-\alpha$. Then $(b-\beta)I=0$. Thus our
 generalized identity (\ref{ps1}) reduces to
 \bea\lb{ps2} ([x_0x_1, x_0x_2]^{2n})(\alpha+\beta)^{n-1}\{1-(\alpha+\beta)^{n}\}(b+\alpha)=0.\eea
 Since this is a trivial (GPI) for $R$, we
  conclude that either $\alpha+\beta=0$ or $(\alpha+\beta)^n=1$ or $b=-\alpha\in C$.

\begin{lem}\label{c3}
Let $R$ be a prime ring with extended centroid $C$ and $I$ be a
right ideal of $R$. Let $H$ be an inner generalized derivation of
$R$. If $H([x,y]^2)^n-H([x,y])^{2n}=0$ for all $x,y\in I$, then
one of the following holds:
\begin{enumerate}
   \item $[I, I]I=0$;
  \item there exists $a\in U$ such that $H(x)=xa$ for all $x\in I$ with $aI=0$;
   \item there exists $a\in U$ such that $H(x)=ax$ for all $x\in R$ with $aI=0$;
  \item there exists $a,b\in U$ such that $H(x)=ax+xb$ for all $x\in R$ with $(a-\alpha)I=(b-\beta)I=0$ for some $\alpha, \beta \in
  C$ and $(\alpha+\beta)^n=1$.
 \end{enumerate}
\end{lem}
%%%%%%%%%%%%%%%%%%%%PROOF%%%%%%%%%%%%%%%%%%%%%%%%%%%%%%%%%%%%%%%%%%%

{\it Proof.} Since $H$ is inner, there exist $a, b\in U$ such that
$H(x)=ax+xb$ for all $x\in R$. If $R$ does not satisfy any
non-trivial (GPI), then by Lemma~\ref{c1}, we conclude that there
exist $\alpha, \beta\in C$ such that $(a-\alpha)I=0$,
$(b-\beta)I=0$ with $\alpha+\beta=0$ or $(\alpha+\beta)^n=1$ or
$b=-\alpha\in C$. If $\alpha+\beta=0$, then for all $x\in I$,
$H(x)=ax+xb=\alpha x+xb=x(\alpha+b)$ with
$0=(\alpha+\beta)I=(\alpha+b)I$, which is our conclusion (2). If
$b=-\alpha\in C$, then for all $x\in R$, $H(x)=ax+xb=(a-\alpha)x$
with $(a-\alpha)I=0$, which is our conclusion (3). In other case
we get our conclusion (4).

 So we assume that $R$ satisfies a non-trivial (GPI).

If $I=R$, then by Lemma \ref{b2}, $a,b\in C$ with $a+b=0$ or
$(a+b)^n=1$. Hence $H(x)=\lambda x$ for all $x\in R$, with
$\lambda^n=1$, since $H$ is nonzero generalized derivation of $R$,
where $\lambda=a+b$. Thus conclusion (4) is obtained.

Now let $I\neq R$. In this case we want to prove that either
$[I,I]I=0$ or there exist $\alpha, \beta\in C$ such that
$(a-\alpha)I=0$ and $(b-\beta)I=0$. To prove this, by
contradiction, we suppose that there exist $c_1, c_2,\cdots,
c_5\in I$ such that
\begin{itemize}
\item $[c_1, c_2]c_3\neq 0;$
 \item $(a-\alpha)c_4\neq 0$ for all $\alpha\in C$ or $(b-\beta)c_5\neq
 0$ for all $\beta\in C$.
 \end{itemize}
 Now we show that this assumption leads a number of
 contradictions.
 Since $R$ satisfies nontrivial (GPI), by \cite{a6},
 $RC$ is a primitive ring having a nonzero socle $H'$ with a nonzero
 right ideal $J=IH'$. Notice that $H'$ is simple, $J=JH'$ and $J$
 satisfies the same basic conditions as $I$. Thus we replace $R$
 by $H'$ and $I$ by $J$.

 Then since $R$ is a regular ring, for $c_1, c_2,\cdots, c_5\in I$
 there exists $e^2=e\in R$ such that
$$eR=c_1R+c_2R+c_3R+c_4R+c_5R.$$
Then $e\in I$ and $ec_i=c_i$  for $i=1, \cdots, 5$. Let $x \in R$.
Then by our hypothesis we have
\begin{equation}\label{c4} (a[e, ex(1-e)]^2+[e,
ex(1-e)]^2b)^n-(a[e, ex(1-e)]+[e, ex(1-e)]b)^{2n}=0.
\end{equation}
Left multiplying by $(1-e)$ we have $((1-e)aex)^{2n}(1-e)=0,$ that
is $((1-e)aex)^{2n+1}=0$ for all $x\in R$. By Levitzki�s lemma
\cite[Lemma 1.1]{her2}, we have $(1-e)aeR=0$ implying $(1-e)ae=0$.
Analogously, right multiplying by $e$, we get $(1-e)be=0$.
Therefore $ae=eae$ and $be=ebe$. Moreover, since $R$ satisfies
 $$e\{(a[ex_1e, ex_2e]^2+[ex_1e, ex_2e]^2b)^n-(a[ex_1e, ex_2e]+[ex_1e,
 ex_2e]b)^{2n}\}e=0,$$ $eRe$ satisfies
 $$(eae[x_1, x_2]^2+[x_1, x_2]^2ebe)^n-(eae[x_1, x_2]+[x_1,
 x_2]ebe)^{2n}=0.$$ Then by Lemma \ref{b2},
one of the following holds: (1) $[eRe, eRe]=0$, (2) $eae, ebe \in
Ce$. Now $[eRe, eRe]=0$ implies $[eR,eR]eR=0$ which contradicts
with the choices of $c_1, c_2, c_3$. Thus $eae=ae\in Ce$ and
$ebe=be\in Ce$. Therefore, there exist $\alpha, \beta\in C$ such
that $(a-\alpha)e=0$ and $(b-\beta)e=0$. This gives
$(a-\alpha)eR=0$ and $(b-\beta)eR=0$. In any case this contradicts
with the choices of $c_4$ and $c_5$.

In case $[I,I]I=0$, conclusion (1) is obtained. Let
$(a-\alpha)I=0$ and $(b-\beta)I=0$ for some $\alpha,\beta\in C$.
Then our hypothesis
 $(a[x,y]^2+[x,y]^2b)^n-(a[x,y]+[x,y]b)^{2n}=0$ for all $x,y\in I$
 gives
 $(\alpha [x,y]^2+[x,y]^2b)^n-(\alpha [x,y]+[x,y]b)^{2n}=0$ for all $x,y\in
 I$. Right multiplying above relation by $[x,y]$, we have
  $(\alpha +\beta)^n\{1-(\alpha +\beta)^n\}[x,y]^{2n+1}=0$ for all $x,y\in
 I$. This implies either $\alpha +\beta=0$ or $(\alpha
 +\beta)^n=1$ or $[x,y]^{2n+1}=0$ for all $x,y\in I$. The last
 relation implies $[I,I]I=0$ (see \cite[Lemma 2 (II)]{chang}), which is our conclusion (1).
 In case $\alpha +\beta=0$, as before, conclusion (2) is obtained. In
 other case conclusion (4) is obtained.

 \vspace*{.4cm}

%%%%%%%%%%%%%%%%%%%%%%%%%%%%%%  THEOREM  %%%%%%%%%%%%%%%%%%%%%%%%%%%%%%%%%%%%%%%%%%%%
Now we are in a position to prove our main theorem for right
ideals.

\vspace*{.2cm}

 {\bf Proof of Theorem \ref{right}.} If $H$ is inner generalized derivation of $R$, then
by Lemma \ref{c3}, we are done. Now let $H$ be not inner. By
Remark 2, we have $H(x)=ax+d(x)$ for some $a\in U$ and a
derivation $d$ on $U$. Let $x,y\in I$. Then by \cite{a1}, $U$
  satisfies
  $$(a[xX,yY]^2+d([xX,yY]^2))^n-(a[xX,yY]+d([xX,yY]))^{2n}=0$$
  that is
  $$(a[xX,yY]^2+d([xX,yY])[xX,yY]+[xX,yY]d([xX,yY]))^n-(a[xX,yY]+d([xX,yY]))^{2n}=0.$$
  This gives  \bea (a[xX,yY]^2+([d(x)X+xd(X),yY]+[xX,d(y)Y+yd(Y)])[xX,yY]\nn\\
  +[xX,yY]([d(x)X+xd(X),yY]+[xX,d(y)Y+yd(Y)]))^n\nn\\
  -(a[xX,yY]+[d(x)X+xd(X),yY]+[xX,d(y)Y+yd(Y)])^{2n}=0.\eea
Since $H$ is not inner, $d$ is also not inner derivation. Then by
Kharchenko's Theorem \cite{a3}, $U$ satisfies
 \bea (a[xX,yY]^2+([d(x)X+xZ_1,yY]+[xX,d(y)Y+yZ_2])[xX,yY]\nn\\
 +[xX,yY]([d(x)X+xZ_1,yY]+[xX,d(y)Y+yZ_2]))^n\nn\\
  -(a[xX,yY]+[d(x)X+xZ_1,yY]+[xX,d(y)Y+yZ_2])^{2n}=0.\eea In
  particular for $X=0$, we have $[xZ_1,yY]^{2n}=0$ for all $Z_1,
  Y\in U.$ In particular,
  $[x,y]^{2n}=0$ for all $x,y\in I$. Then by \cite[Lemma 2 (II)]{chang}, $[I,I]I=0$,
  which is our conclusion (1).\\

  From above Theorem \ref{right} following
corollaries  are straightforward.
 \begin{cor}\label{gsg}
Let $R$ be a prime ring, $I$ a non-zero right ideal of $R$ and $H$
a non-zero generalized derivation of $R$. If $H$ acts as a Jordan
homomorphism on the set $[I, I]$, then one of the following holds:
\begin{enumerate}
  \item $[I, I]I=0$;
  \item there exists $a\in U$ such that $H(x)=xa$ for all $x\in I$ with $aI=0$;
  \item there exists $a\in U$ such that $H(x)=ax$ for all $x\in R$ with $aI=0$;
  \item there exists $q\in U$ such that $H(x)=xq$ for all $x\in I$ with $qx=x$ for all $x\in
  I$.
\end{enumerate}
\end{cor}
 {\it Proof.} By Theorem \ref{right}, conclusions (1)-(3) are
 obtained. Thus we have only to consider the case, when $H(x)=ax+xb$ for all $x\in R$ with $(a-\alpha)I=(b-\beta)I=0$ for some $\alpha, \beta \in
  C$ and $\alpha+\beta=1$. In this case, for all $x\in I$, we have
  $H(x)=ax+xb=\alpha x+xb=x(\alpha+b)$, where
  $0=(b-\beta)I=(b+\alpha-1)I$. This is our conclusion (4).

  \singlespacing

\vspace*{.4cm}

\nii {\small Basudeb Dhara \\
 Department of Mathematics\\  Belda College, Belda\\
 Paschim Medinipur-721424, INDIA\\  e-mail:
basu$\_$dhara@yahoo.com\\ \\

\nii Shervin Sahebi, Venus Rahmani\\
Department of Mathematics\\ Islamic Azad University\\
 Central Tehran Branch, 13185/768, Tehran, IRAN\\
 e-mail: sahebi@iauctb.ac.ir\\
 e-mail: ven.rahmani.math@iauctb.ac.ir}


\small


\begin{thebibliography}{1}

\bibitem{a01}
A. Asma, N. Rehman, A. Shakir, On Lie ideals with derivations as
homomorphisms and anti-homomorphisms, \emph{ Acta Math. Hungar.}
101 (1-2) (2003), 79-82.


\bibitem{a0}
 K. I. Beidar, W. S.  Martindale III, A. V. Mikhalev,   \emph{Rings with generalized identities}.
 Pure and Applied Math. Vol. 196 (1996), New York, Marcel Dekker.

 \bibitem{a010}
 J. Bergen, I. N. Herstein, J. W. Kerr,  Lie ideals and derivations of prime rings, \emph{J. Algebra.} 71 (1981),
 259-267.

  \bibitem{a1}
 C. L. Chuang,  GPI's having coefficients in Utumi quotient rings
 \emph{proc. Amer. Math. soc.} 103 (1988), 723-728.


\bibitem{chang}
C. M. Chang,  Power central values of derivations on multilinear
polynomials, \emph{Taiwanese J. Math.}, 7 (2) (2003), 329-338.

\bibitem {cl}
C. M. Chang and T. K. Lee, Annihilators of power values of
derivations in prime rings,  \emph{Comm. Algebra. }, 26 (7)
(1998), 2091-2113.

\bibitem {filippis} V. De Filippis, Generalized derivations as Jordan homomorphisms
 on Lie ideals and right ideals, {\it Acta Math. Sinica, English
Series} 25 (2) (2009), 1965-1974.

\bibitem{u54}
T. S. Erickson,  W. S. Martindale III, J. M. Osborn,  Prime
nonassociative  algebras, \emph{ Pacific J.
 Math.} 60 (1975), 49-63.

 \bibitem {x11}
 C. Faith, Y. Utumi, On a new proof of Litoff's
theorem, \emph{ Acta Math. Acad. Sci. Hung.} 14 (1963), 369-371.


\bibitem {Gol}
O. Golbasi, K. Kaya, On Lie ideals with generalized derivations,
\emph{Sibrian Math.}, 47 (5) (2006), 862-866.


   \bibitem{her2}
 I. N. Herstein,  \emph{Topics in ring theory}.
   Univ. of Chicago Press, Chicago, (1969).

 \bibitem{a5}
N. Jacobson,  \emph{Structure of rings}.
 Amer. Math. Soc. Colloq. Pub. 37. Providence, RI: Amer. Math. Soc.
 (1964).

 \bibitem{a3}
 V. K. Kharchenko, Differential identity of prime rings,
 \emph{Algebra and Logic} 17 (1978), 155-168.

  \bibitem{a03}
 C. Lanski,  An engle condition with derivation,
 \emph{Proc. Amer. Math. Soc.} 183 (3) (1993), 731-734.

\bibitem{lm}
C. Lanski,  S. Montgomery, Lie structure of prime rings of
characteristic $2$,  \emph{Pacific J. Math.} 42 (1) (1972),
117-136.

\bibitem{g}
T. K. Lee, Semiprime rings with differential identities,
\emph{Bull. Inst. Math. Acad. Sinica}, 20 (1) (1992), 27-38.

 \bibitem{a6}
 W. S. Martindale III,  Prime rings satistying a generalized polynomial
  identity, \emph{J. Algebra.} 12 (1972), 576-584.



\end{thebibliography}
\end{document}